\documentclass[12pt,reqno]{amsart}

\usepackage{latexsym}
\usepackage{amssymb}
\usepackage{amsmath}
\usepackage{fancybox,color}

 \def\1{\raisebox{2pt}{\rm{$\chi$}}}

\newtheorem{theorem}{Theorem}[section]

\newtheorem{lemma}[theorem]{Lemma}
\newtheorem{proposition}[theorem]{Proposition}
\newtheorem{definition}[theorem]{Definition}
\newtheorem{remark}[theorem]{Remark}

\newcommand{\RR}{{\mathbb R}}

 \def\1{\raisebox{2pt}{\rm{$\chi$}}}

\usepackage{enumerate}

\renewcommand{\div}{\operatorname{div}}

\begin{document}

\dedicatory{To the memory of Oded Schramm}
\title[$p(x)$-harmonic functions]
{$p(x)$-harmonic functions with unbounded exponent\\
in a subdomain}

\author{J.~J.~Manfredi, J.~D.~Rossi and J.M.~Urbano}

\address{Juan J.~Manfredi
\hfill\break\indent Department of Mathematics
\hfill\break\indent University of Pittsburgh \hfill\break\indent
Pittsburgh PA 15260, USA
\hfill\break\indent {\tt
manfredi@pitt.edu}}
\address{Julio D.~Rossi
\hfill\break\indent IMDEA Matem{\'a}ticas, \hfill\break\indent C-IX,
Campus UAM, \hfill\break\indent Madrid, Spain \hfill\break\indent On
leave from Departamento de Matem\'{a}tica, \hfill\break\indent FCEyN
UBA (1428) \hfill\break\indent Buenos Aires, Argentina.
\hfill\break\indent {\tt jrossi@dm.uba.ar}}
\address{Jos{\'e} Miguel Urbano
\hfill\break\indent CMUC, Department of Mathematics
\hfill\break\indent University of Coimbra \hfill\break\indent
3001-454 Coimbra, Portugal \hfill\break\indent {\tt
jmurb@mat.uc.pt}}

\keywords{$p(x)$-Laplacian, infinity-Laplacian, viscosity solutions.\\
\indent 2000 {\it Mathematics Subject Classification.} 35J20, 35J60,
35J70.}

\begin{abstract}
We study the Dirichlet problem $-\div(|\nabla u|^{p(x)-2} \nabla u)
=0 $ in $\Omega$, with $u=f$ on $\partial \Omega$ and $p(x) =
\infty$ in $D$, a subdomain of the reference domain $\Omega$. The
main issue is to give a proper sense to what a solution is. To this
end, we consider the limit as $n \to \infty$ of the solutions $u_n$
to the corresponding problem when $p_n(x) =p(x) \wedge n$, in
particular, with $p_n = n$ in $D$. Under suitable assumptions on the
data, we find that such a limit exists and that it can be
characterized as the unique solution of a variational minimization
problem which is, in addition, $\infty$-harmonic within $D$.
Moreover, we examine this limit in the viscosity sense and find the
boundary value problem it satisfies in the whole of $\Omega$.
\end{abstract}

\maketitle

\section{Introduction}

The goal of this paper is to study the elliptic problem
\begin{equation} \label{1.1}
\left\{\begin{array}{ll} \displaystyle -\Delta_{p(x)} u (x) =0,
\qquad & x\in\Omega\subset\mathbb{R}^{N},\\
[5pt]
 u(x) = f(x),\quad & x\in \partial \Omega,
\end{array}\right.
\end{equation}
where $\Delta_{p(x)} u (x) := \mbox{div} \left( |\nabla u (x)|^{p(x)
-2} \nabla u(x) \right)$ is the $p(x)$-Laplacian operator and the
variable exponent $p(x)$ verifies
\begin{equation}\label{p=infty}
    p(x) = +\infty, \qquad x \in D ,
    \end{equation}
for some subdomain $D\subset \Omega$. We assume that $\Omega$ and
$D$ are bounded and convex domains with smooth boundaries, at
least of class $C^1$. On the complementary domain $\Omega
\setminus
\overline D$ we assume that $p(x)$ is a continuously
differentiable bounded function.

On the variable exponent, apart from \eqref{p=infty}, we also
require that
\begin{equation}\label{p>N}
p_- := \inf_{x \in \Omega} p(x) > N, \end{equation} so that we will
always be dealing with continuous solutions for \eqref{1.1}; to fix
notation, we define
$$p_+ := \sup_{x \in \Omega \setminus \overline D} p(x).$$

The boundary data $f$ is taken to be Lipschitz continuous.

Our strategy to solve \eqref{1.1} is to replace $p(x)$ by a sequence
of bounded functions $p_n(x)$ such that $p_n (x)$ is increasing and
converging to $p(x)$.  For definiteness, we consider, for $n>N$,
$$
p_n (x) := \min\{p(x), n\}.
$$
We will use the notation \eqref{1.1}$_n$ to refer to problem
\eqref{1.1} for the variable exponents $p_n(x)$.

Since $p(x)$ is bounded in $\Omega \setminus D$, we have, for
large $n$, specifically for $n> p_+$,
$$p_n (x)= \left\{
\begin{array}{ll}
p(x), \quad & x \in \Omega \setminus D, \\[5pt]
n, \quad & x \in D.
\end{array} \right.
$$
Moreover, still for large $n$, the boundary of the set $\{ p(x)>n
\}$ coincides with the boundary of $D$ and thus does not depend on
$n$. This fact is important when passing to the limit.

Using a variational method, we solve \eqref{1.1}$_n$ obtaining
solutions $u_n$; if the limit
\begin{equation}\label{limite.n}
    \lim_{n \to \infty } u_n
\end{equation}
exists, we call it $u_\infty$. It is a natural candidate to be a
solution to \eqref{1.1} with the original variable exponent $p(x)$.
A crucial role in this process will be played by the set
$$
S = \left\{ u \in W^{1,p^{-}}(\Omega) \ : \ u|_{\Omega \setminus
\overline D} \in W^{1, p(x)} (\Omega \setminus \overline D), \quad
\| \nabla u \|_{L^\infty (D)} \leq 1 \quad \mbox{and} \quad u
|_{\partial \Omega}= f \right\}
$$
and by  the infinity Laplacian
$$\Delta_{\infty} u : =  \left( D^2 u \: \nabla u \right) \cdot \nabla u = \sum_{i,j=1}^N
\frac{\partial u}{\partial x_i} \frac{\partial u}{\partial x_j}
\frac{\partial^2 u}{\partial x_i x_j}.$$
Our main results are condensed in the following theorem.

\medskip

\noindent \textbf{Theorem.} \textit{There exists a unique solution
$u_n$ to \eqref{1.1}$_n$. If $S\neq \emptyset$, then the uniform
limit
$$u_\infty := \lim_{n\to \infty} u_n$$
exists and is characterized as the unique function that is a
minimizer of the integral
\begin{equation}\label{referee1}
\displaystyle \int_{\Omega \setminus \overline D} \frac{|\nabla
u|^{p(x)} }{p(x)} \, dx
\end{equation}
in $S$ and, in addition, verifies $$-\Delta_\infty u_\infty =0 \quad
\mbox{ in } D,$$ in the viscosity sense. Moreover, $u_\infty$ is a
viscosity solution of
$$
\left\{\begin{array}{ll} - \Delta_{p(x)} u (x) =0, \qquad &
x\in\Omega
\setminus \overline D,\\
[5pt]- \Delta_{\infty} u (x) =0, \qquad & x\in D,\\
[5pt]  \mathrm{sgn} (|\nabla u | (x) -1) \: \mathrm{sgn}  \left(
\frac{\partial u}{\partial \nu}(x) \right) =0, \qquad & x \in
\partial D
\cap \Omega, \\
[5pt] u(x) = f(x),\quad & x\in \partial \Omega,
\end{array}\right.
$$
where $\nu$ is the exterior unit normal vector to $\partial D$ in
$\Omega$. }\par
\textit{
Finally, if $\partial \Omega \cap \overline{D} \neq \emptyset$ and
the Lipschitz constant of $f|_{\partial \Omega \cap \overline{D}}$
is strictly greater than one, then $S = \emptyset$ and we have
$$
\liminf_{n \to \infty} \left( \int_D \frac{|\nabla u_n|^{n} }{n} \, dx
\right)^{1/n}>1;
$$
hence, the natural energy associated to $u_n$ is unbounded.}

\medskip

\begin{remark}
The boundedness of $\Omega$ is used to ensure compactness of
minimizing sequences for {\rm (\ref{referee1})}, while the
convexity of $\Omega$ and $D$ guarantees that the Lipschitz
constant of $W^{1,\infty}$ functions coincides with the
$L^{\infty}$ norm of their gradients, which will be instrumental
in some of the proofs.
\end{remark}

\begin{remark}
The characterization of the non-emptiness of $S$ is an interesting
open problem that strongly depends on the geometry of $\Omega$ and
$D$, and on the boundary data $f$. When $\partial \Omega \cap
\overline{D} = \emptyset$, $S$ is always non-empty. When $\partial
\Omega \cap \overline{D} \neq \emptyset$, the condition that the
Lipschitz constant of $f|_{\partial \Omega \cap \overline{D}}$ is
less than or equal to one is necessary but, in general, it is not
sufficient (cf. section \ref{remark222}).
\end{remark}

Partial differential equations involving variable exponents became
popular a few years ago in relation to applications to elasticity
and electrorheological fluids. Meanwhile, the underlying functional
analytical tools have been extensively developed and new
applications, \textit{e.g.} to image processing, have kept the
subject as the focus of an intensive research activity. For general
references on the $p(x)$-Laplacian we refer to \cite{DHN}, that
includes a thorough bibliography, and \cite{KR}, a seminal paper
where many of the basic properties of variable exponent spaces were
established. The delicate regularity properties of $p(x)$-harmonic
functions have been established in \cite{AcMi1} and \cite{AcMi2}.

In the literature, the variable exponent $p(x)$ is always assumed
to be bounded, a necessary condition to define a proper norm in
the corresponding Lebesgue spaces. To the best of our knowledge,
this paper is the first attempt at analyzing a problem where the
exponent $p(\cdot)$ becomes infinity in some part of the domain.
For constant exponents, limits as $p\to \infty$ in $p-$Laplacian
type problems have been widely studied, see for example
\cite{BBM}, and are related to optimal transport problems
(\textit{cf.} \cite{Am}).

\medskip

{\bf Organization of the paper.} The rest of the paper is organized
as follows: in Section \ref{NLPpg1} we show existence and uniqueness
of solutions with $p(x) =p_n (x) = p \wedge n$ using a variational
argument; moreover we find the equation that they verify in the
viscosity sense and prove some useful independent of $n$ estimates;
in Section \ref{sect-variational} we pass to the limit in the
variational formulation of the problem and we deal with the limit in
the viscosity sense; in Section \ref{remark222} we discuss necessary and sufficient conditions related to the non-emptiness of $S$ and present examples and counter-examples. Finally, in Section \ref{1-d} we present a
detailed analysis of the one-dimensional case.

\section{Weak and viscosity approximate solutions}\label{NLPpg1}
\setcounter{equation}{0}

To start with, let us establish the existence and uniqueness of the
approximations $u_n$ in the weak sense.

\begin{lemma} \label{existencia.un}
There exists a unique weak solution $u_n$ to \eqref{1.1}$_n$,
which is the unique minimizer of the functional
\begin{equation}\label{funcional.p(x).2}
    F_n(u) = \int_\Omega \frac{|\nabla u|^{p_n (x)} }{p_n (x)} \, dx
    = \int_D \frac{|\nabla u|^{n } }{n} \, dx + \int_{\Omega \setminus \overline D}
    \frac{|\nabla u|^{p (x)} }{p (x)} \, dx
\end{equation}
in
\begin{equation}\label{Sn.2}
    S_n = \left\{ u \in W^{1, p_n (x)} (\Omega) \, : \, u |_{\partial
\Omega}= f \right\}.
\end{equation}
\end{lemma}

\begin{proof}
Although the exponent $p_n (\cdot)$ might be discontinuous,
functions in the variable exponent Sobolev space $W^{1, p_n (\cdot)}
(\Omega)$ are continuous thanks to assumption \eqref{p>N}. Indeed,
for $n$ sufficiently large, we have $p_n (\cdot)\geq (p_n)_- \geq
p_-
> N$ and the continuous embedding in
\begin{equation}\label{embedd}
W^{1, p_n (\cdot)} (\Omega) \hookrightarrow W^{1, p_-} (\Omega)
\subset C \left( \overline{\Omega} \right)
\end{equation}
follows from \cite[Theorem 2.8 and (3.2)]{KR}. That the boundedness
away from the dimension is not superfluous when the exponent is not
continuous is shown by a counter-example in \cite[Example 3.3]{HH1}.

We can then take the boundary condition $u |_{\partial \Omega}= f$
in the classical sense (recall that $f$ is assumed to be Lipschitz)
and the results of \cite{HHKV} apply since the jump condition (cf.
\cite[(4.1)-(4.2)]{HHKV}) is trivially satisfied by the variable
exponent because $p_n (\cdot)\geq N$. This is a sufficient condition
for a $p_n (\cdot)$-Poincar{\'e} inequality to hold in $W_0^{1, p_n
(\cdot)} (\Omega)$ which, in turn, is instrumental in obtaining the
coercivity of the functional. The lower semicontinuity is standard
as is the strict convexity, that also gives the uniqueness.

It is also standard that the minimizer of $F_n$ in $S_n$ is the
unique weak solution of \eqref{1.1}$_n$, \textit{i.e.}, $u_{n}=f$
on $\partial \Omega$ and it satisfies the weak form of the
equation, namely,
\begin{equation}\label{ec.debil.1}
\int_\Omega  |\nabla u_n |^{p_n(x) -2} \nabla u_n \cdot \nabla
\varphi \, dx=0, \qquad \forall \: \varphi \in C^\infty_0 (\Omega).
\end{equation}
\end{proof}

\begin{lemma} \label{lema.sol.debil}
Problem \eqref{1.1}$_n$ can be rewritten as
\begin{equation} \label{1.1.n.2}
\left\{\begin{array}{ll} - \Delta_{p(x)} u_n (x) =0,\qquad &
x\in\Omega \setminus \overline D,\\
[5pt]- \Delta_{n} u_n (x) =0, \qquad & x\in D, \\
[5pt]
 \displaystyle |\nabla u_n (x)|^{n -2} \frac{\partial u_n }{\partial \nu}(x)
 = |\nabla u_n (x)|^{p (x) -2} \frac{\partial u_n }{\partial
\nu}(x),
\qquad & x\in \partial D \cap \Omega,\\
[5pt]
u_n(x) = f(x),\quad & x\in \partial \Omega,
\end{array}\right.
\end{equation}
where $\nu$ is the exterior unit normal to $\partial D$ in $\Omega$.
\end{lemma}

\begin{proof}
Just notice that the weak form of this problem is exactly the same
as the one that holds for \eqref{1.1}$_n$. This follows since
after multiplying by a test function and integrating by parts one
arrives at
\eqref{ec.debil.1} for both problems.
\end{proof}

Next, we investigate the problem satisfied by $u_n$ from the point
of view of viscosity solutions.

Let us recall the definition of viscosity solution (see \cite{CIL}
and \cite{Bar}) for a problem like \eqref{1.1.n.2}, which involves a
transmission condition across the boundary $\partial D \cap \Omega$.
Assume we are given a family of continuous functions
$$F_i:\overline{\Omega}\times\mathbb{R}^N\times\mathbb{S}^{N\times N}\rightarrow\mathbb{R}.$$
 The associated equations
$$F_i(x, \nabla u, D^2 u)=0$$
are called (degenerate) elliptic if
$$F_i(x,\xi, X)\le F_{i}(x,\xi, Y)\quad \textrm{whenever} \ X\ge Y.$$

\begin{definition} \label{def.sol.viscosa}
Consider the problem
\begin{equation}\label{ec.viscosa.con.borde}
\begin{array}{ll}
F_1 (x, \nabla u, D^2 u )   = 0, \qquad & \textrm{in} \ \: \Omega \setminus \overline D, \\
[5pt] F_2 (x, \nabla u, D^2 u )   = 0, \qquad & \textrm{in} \ \:
D,
\end{array}
\end{equation}
with a transmission condition
\begin{equation}\label{ec.transmission}
\begin{array}{ll}
B(x, u , \nabla u )  = 0, \qquad & \textrm{on} \ \: \partial D \cap
\Omega,
\end{array}
\end{equation}
and a boundary condition
\begin{equation}\label{boundary.cond}
\begin{array}{ll}
u  = f, \qquad & \textrm{on} \ \: \partial \Omega.
\end{array}
\end{equation}
A lower semi-continuous function $u$ is a viscosity supersolution of
\eqref{ec.viscosa.con.borde}--\eqref{boundary.cond} if $u\geq f$ on
$\partial \Omega$ and for every $\phi \in C^2(\overline{\Omega})$
such that $u-\phi $ has a strict minimum at the point $x_0 \in
{\Omega}$, with $u(x_0)= \phi(x_0)$, we have
$$\begin{array}{rcccl} F_1(x_0, \nabla \phi (x_0), D^2\phi (x_0)) &
\ge & 0 & \textrm{if} &
x_0 \in \Omega \setminus \overline D,\\
[10pt] F_2(x_0, \nabla \phi (x_0), D^2\phi (x_0)) & \ge & 0 &
\textrm{if} & x_0 \in D,\\
[10pt] \max \left\{
\begin{array}{l}
F_1 \left( x_0, \nabla \phi (x_0), D^2\phi (x_0) \right) \\
[5pt]
F_2(x_0, \nabla \phi (x_0), D^2\phi (x_0)) \\
[5pt]
B (x_0, \phi (x_0), \nabla \phi (x_0))
\end{array}
\right\} & \ge & 0 & \textrm{if} & x_0 \in \partial D \cap \Omega.
\end{array}$$
An upper semi-continuous function $u$ is a viscosity subsolution of
\eqref{ec.viscosa.con.borde}--\eqref{boundary.cond} if $u\le f$ on
$\partial \Omega$ and for every  $\psi \in C^2(\overline{\Omega})$
such that $u-\psi $ has a strict maximum at the point $x_0 \in
{\Omega}$, with $u(x_0)= \psi(x_0)$, we have
$$\begin{array}{rcccl}
F_1(x_0, \nabla \psi (x_0), D^2\psi (x_0)) & \le & 0 & \textrm{if} &
x_0 \in \Omega \setminus \overline D,\\
[10pt] F_2(x_0, \nabla \psi (x_0), D^2\psi (x_0)) & \le & 0 &
\textrm{if} & x_0 \in D,\\
[10pt] \min \left\{
\begin{array}{l}
F_1 \left( x_0, \nabla \psi (x_0), D^2\psi (x_0) \right) \\
[5pt]
F_2(x_0, \nabla \psi (x_0), D^2\psi (x_0)) \\
[5pt] B (x_0, \psi (x_0), \nabla \psi (x_0))
\end{array}
\right\} & \le & 0 & \textrm{if} & x_0 \in \partial D \cap \Omega.
\end{array}$$
Finally, $u$ is a viscosity solution if it is both a viscosity
supersolution and a viscosity subsolution.
\end{definition}

In the sequel, we will use the notation as in the definition: $\phi$
will always stand for a test function touching the graph of $u$ from
below and $\psi $ for a test function touching the graph of $u$ from
above.

\begin{proposition} \label{lema.sol.debil.es.sol.viscosa}
Let $u_n$ be a continuous weak solution of \eqref{1.1}$_n$. Then
$u_n$ is a viscosity solution of \eqref{1.1.n.2} in the sense of
Definition \ref{def.sol.viscosa}.
\end{proposition}

\begin{proof}
To simplify, we omit in the proof the subscript $n$. Let $x_0 \in
\Omega\setminus \overline D$ and a let $\phi$ be a test function
such that $u(x_0)=\phi (x_0)$ and $u-\phi$ has a strict minimum at
$x_0$. We want to show that
\begin{eqnarray*}
- \Delta_{p(x_0)} \phi (x_0) & = & -|\nabla
\phi(x_{0})|^{p(x_{0})-2} \Delta \phi(x_{0})-(p(x_{0})-2)
|\nabla \phi(x_{0})|^{p(x_{0})-4} \Delta_\infty \phi(x_{0}) \\
& & -|\nabla \phi(x_{0})|^{p(x_{0})-2} \ln (|\nabla \phi|)(x_{0})
\left\langle \nabla \phi(x_{0}) , \nabla p(x_{0}) \right\rangle
\\
&\ge & 0.
\end{eqnarray*}

Assume, \textit{ad contrarium},  that this is not the case; then
there exists a radius $r>0$ such that $B(x_0,r) \subset
\Omega\setminus \overline D$ and
\begin{eqnarray*}
- \Delta_{p(x)} \phi (x) & = & -|\nabla \phi(x)|^{p(x)-2} \Delta \phi(x)-(p(x)-2)|\nabla \phi(x)|^{p(x)-4} \Delta_\infty \phi(x) \\
& & -|\nabla \phi(x)|^{p(x)-2} \ln (|\nabla \phi|)(x) \langle \nabla
\phi(x) , \nabla p(x) \rangle
\\
& < & 0,
\end{eqnarray*}
for every $x\in B(x_0,r)$. Set $m = \inf_{|x-x_0|=r} (u-\phi)(x)$
and let $\Phi (x) = \phi(x) + m/2$. This function $\Phi$ verifies
$\Phi (x_0) > u(x_0)$ and
\begin{equation}
- \Delta_{p(x)} \Phi = -\mbox{div} ( |\nabla \Phi|^{p(x)-2} \nabla
\Phi) <0 \qquad \textrm{in} \ \: B(x_{0},r). \label{p-lap}
\end{equation}

\noindent Multiplying \eqref{p-lap} by $(\Phi - u)^+$, which
vanishes on the boundary of $B(x_0,r)$, we get
$$
\int_{B(x_{0},r) \cap \{ \Phi > u \}} |\nabla \Phi|^{p(x)-2} \nabla
\Phi \cdot \nabla (\Phi - u) \, dx < 0.
$$
On the other hand, taking $(\Phi - u)^+$, extended by zero
outside $B(x_{0},r)$, as test function in the weak formulation of
\eqref{1.1}$_n$, we obtain
$$
\int_{ B(x_{0},r) \cap \{ \Phi > u \} } |\nabla u|^{p(x)-2} \nabla u
\cdot \nabla (\Phi - u) \, dx=0,
$$
since $p_n(x) = p(x)$ in $\Omega\setminus \overline D$. Upon
subtraction and using a well know inequality, see for example \cite{L}, we
conclude
\begin{eqnarray*}
0 & > & \int_{  B(x_{0},r) \cap \{ \Phi
> u \}} \left( |\nabla \Phi|^{p(x)-2} \nabla \Phi - |\nabla u|^{p(x)-2}
\nabla u \right) \cdot \nabla (\Phi - u) \, dx\\
& \geq & c \int_{ B(x_{0},r) \cap \{ \Phi > u \}} |\nabla \Phi -
\nabla u|^{p(x)} \, dx,
\end{eqnarray*}
a contradiction. Here $c$ is a constant that depends on $N$, $p^{-}$
and $\displaystyle \sup_{x\in B(x_{0},r)} p(x)$.

If $x_0 \in D$ the proof is entirely analogous, albeit simpler due
to the absence of the logarithmic term, and we obtain
$$- \Delta_n
\phi (x_0)  =  -|\nabla \phi(x_{0})|^{n-2} \Delta \phi(x_{0})-(n-2)
|\nabla \phi(x_{0})|^{n-4} \Delta_\infty \phi(x_{0}) \ge  0.
$$
The constant $c$ in this case depends on $N$ and $n$.

If $x_0 \in \partial D \cap \Omega$ we want to prove that
$$\max \left\{
\begin{array}{l}
- \Delta_{p(x_0)} \phi (x_0) \\
[5pt] - \Delta_n \phi (x_0) \\
[5pt]   | \nabla \phi (x_0)|^{n-2} \frac{\partial \phi}{\partial
\nu}(x_0) -| \nabla \phi (x_0)|^{p(x_0)-2} \frac{\partial
\phi}{\partial \nu}(x_0)
\end{array}
\right\}  \ge  0.
$$
If this is not the case, there exists a radius $r>0$ such that
$$
- \Delta_{p(x)} \phi (x) < 0 \qquad \textrm{and} \qquad - \Delta_n
\phi (x) <0,
$$
for every $x\in B(x_0,r)$. Set $m = \inf_{|x-x_0|=r} (u-\phi)(x)$
and let $\Phi (x) = \phi(x) + m/2$. This function $\Phi$ verifies
$\Phi (x_0) > u(x_0)$,
\begin{equation} \label{eqq1}
- \Delta_{p(x)} \Phi <0 \qquad \textrm{in} \ \: B(x_0, r) \cap
(\Omega \setminus \overline{D})
\end{equation}
and
\begin{equation} \label{eqq2}
- \Delta_n \Phi <0 \qquad \textrm{in} \ \:  B(x_0, r) \cap D.
\end{equation}
Moreover, we can assume (taking $r$ smaller if necessary) that
\begin{equation}\label{cond.borde.tranferencia}
| \nabla \Phi (x)|^{n-2} \frac{\partial \Phi}{\partial \nu}(x) - |
\nabla \Phi (x)|^{p(x)-2} \frac{\partial \Phi}{\partial \nu}(x) <0
\qquad \textrm{in} \ \:  B(x_0, r) \cap \partial D .
\end{equation}

Multiplying both \eqref{eqq1} and \eqref{eqq2} by $(\Phi - u)^+$,
integrating by parts and adding, we obtain
$$\int_{B(x_0, r) \cap \Omega \setminus \overline{D}} |\nabla
\Phi|^{p(x)-2} \nabla \Phi \cdot \nabla (\Phi - u)^+ \, dx +
\int_{B(x_0, r) \cap D} |\nabla \Phi|^{n-2} \nabla \Phi \cdot \nabla
(\Phi - u)^+ \, dx$$ $$< \int_{B(x_0, r) \cap \partial D} \left( |
\nabla \Phi |^{n-2} \frac{\partial \Phi}{\partial \nu} - | \nabla
\Phi|^{p(x)-2} \frac{\partial \Phi}{\partial \nu} \right) (\Phi -
u)^+ \, dS,
$$
taking also into account that the test function vanishes on the
boundary of $B(x_0,r)$. Using \eqref{cond.borde.tranferencia}, we
finally get
$$\hspace*{-3cm}\int_{B(x_0, r) \cap (\Omega \setminus \overline{D}) \cap \{ \Phi > u \}} |\nabla
\Phi|^{p(x)-2} \nabla \Phi \cdot \nabla (\Phi - u) \, dx $$
$$\hspace*{5cm}+ \int_{B(x_0, r) \cap D \cap \{ \Phi > u \}} |\nabla
\Phi|^{n-2} \nabla \Phi \cdot \nabla (\Phi - u) \, dx < 0.$$ On the
other hand, taking $(\Phi - u)^+$, extended by zero outside
$B(x_{0},r)$, as test function in the weak formulation of
\eqref{1.1}$_n$, we reach a contradiction as in the previous cases.
This proves that $u$ is a viscosity supersolution.

The proof that $u$ is a viscosity subsolution runs as above and we
omit the details.
\end{proof}

We next obtain uniform estimates (independent of $n$) for the
sequence of approximations $(u_n)_n$.

\begin{proposition} \label{lemma.estimaciones}
Assume the set
$$
S = \left\{ u \in W^{1,p^{-}}(\Omega) \ : \ u|_{\Omega \setminus
\overline D} \in W^{1, p(x)} (\Omega \setminus \overline D), \ \ \|
\nabla u \|_{L^\infty (D)} \leq 1 \ \ \mbox{and} \ \ u |_{\partial
\Omega}= f \right\}
$$
is nonempty. Then $u_n$, the minimizer of $F_n$ in $S_n$, satisfies
$$
F_n (u_n)= \int_{\Omega} \frac{|\nabla u_n|^{p_n (x)} }{p_n (x)} \,
dx \leq \int_D \frac{|\nabla v|^{n} }{n} \, dx + \int_{\Omega
\setminus \overline D}\frac{|\nabla v|^{p (x)} }{p (x)} \, dx,
$$
for every $v\in S$. Hence, the sequence $(F_n (u_n))_n$ is uniformly
bounded and the sequence $(u_n)_n$ is uniformly bounded in $W^{1,
p_-} (\Omega)$ and equicontinuous.
\end{proposition}

\begin{proof}  Recalling \eqref{Sn.2}, the definition of $S_n$, observe
that $S \subset S_{n}$, for every $n$. Since $u_n$ is a minimizer,
we have
$$
F_n (u_n) \leq F_n (v), \quad \forall \ v \in S.$$ Hence, picking an
element $v \in S\neq \emptyset$,
\begin{eqnarray*}
F_n(u_n) = \int_\Omega \frac{|\nabla u_n|^{p_n(x)} }{p_n(x)} \, dx
&\leq &\int_\Omega \frac{|\nabla v|^{p_n(x)} }{p_n(x)} \, dx \\
&=&
\int_D \frac{|\nabla v|^{n} }{n} \, dx +\int_{\Omega \setminus
\overline D}\frac{|\nabla v|^{p (x)} }{p (x)} \, dx\\
& \leq & |D| + \int_{\Omega \setminus \overline D}\frac{|\nabla
v|^{p (x)} }{p (x)} \, dx \equiv C_\ast.
\end{eqnarray*}

In order to estimate the Sobolev norm, we first use Poincar{\'e}
inequality and the boundary data, to obtain
\begin{eqnarray*}
\| u_n\|_{W^{1, p_-} (\Omega)} & \leq & \| u_n -f\|_{W_0^{1,
p_-} (\Omega)} + \| f\|_{W^{1, p_-} (\Omega)} \\
& \leq & C \, \|
\nabla (u_n -f)\|_{L^{p_-} (\Omega)} + \| f\|_{W^{1, \infty} (\Omega)}\\
& \leq & C \, \| \nabla u_n\|_{L^{p_-} (\Omega)} + (C+1)   \|
f\|_{W^{1, \infty} (\Omega)}.
\end{eqnarray*}
We proceed, using H\"{o}lder inequality and elementary computations,
to get
\begin{eqnarray*}
\| \nabla u_n\|_{L^{p_-} (\Omega)} & = & \left( \int_\Omega |\nabla
u_n|^{p_-} \, dx \right)^{1/p_-} \leq \left( \int_D |\nabla
u_n|^{p_-} \, dx \right)^{1/p_-} + \left( \int_{\Omega \setminus
\overline D} |\nabla u_n|^{p_-} \, dx \right)^{1/p_-}\\
& = & \left( \int_D |\nabla u_n|^{p_-} \, dx \right)^{1/p_-} +
\left( \int_{(\Omega \setminus \overline D) \cap  \{|\nabla u_n|\leq
1\}} |\nabla u_n|^{p_-} \, dx \right)^{1/p_-}
\\
& & \displaystyle + \left( \int_{(\Omega \setminus
\overline D) \cap  \{|\nabla u_n| > 1\}} |\nabla u_n|^{p_-} \, dx \right)^{1/p_-}\\
& \leq & |D|^{\frac{1}{p_-}-\frac{1}{n}} \left( \int_{D} | \nabla
u_n |^n \right)^{1/n} + |\Omega| + \left( \int_{\Omega \setminus
\overline D} |\nabla u_n|^{p(x)} \, dx \right)^{1/p_-}.
\end{eqnarray*}
Since we have the bounds
$$\left( \int_{D} | \nabla
u_n |^n \right)^{1/n} = n^{1/n} \left( \int_D \frac{|\nabla
u_n|^n}{n} \, dx \right)^{1/n} \leq n^{1/n}\left( F_n(u_n)
\right)^{1/n} \leq 2 C_{\ast}$$
and
$$ \int_{\Omega \setminus
\overline D} |\nabla u_n|^{p(x)} \, dx \leq p_+ \int_{\Omega
\setminus \overline D} \frac{|\nabla u_n|^{p(x)}}{p(x)}  \, dx \leq
p_+ F_n(u_n) \leq p_+ C_{\ast},$$ we conclude that the sequence  $(u_n)_n$ is
uniformly bounded in $W^{1, p_-} (\Omega)$ and, recalling the
embedding in \eqref{embedd}, that it is equicontinuous.
\end{proof}

\section{Variational and viscosity limit}\label{sect-variational}
\setcounter{equation}{0}

We first analyze the case in which $\partial \Omega \cap \overline{D} \neq \emptyset$ and the Lipschitz constant of
$f|_{\partial \Omega \cap \overline{D}}$ is greater than one. Note
that, in this case, $S=\emptyset$ since any Lipschitz
extension $u$ of this datum to $D$ verifies $\| \nabla u
\|_{L^\infty (D)} >1$.

\begin{theorem}\label{no.acotada}
Assume that $\partial \Omega \cap \overline{D} \neq \emptyset$ and the Lipschitz constant of
$f|_{\partial \Omega \cap \overline{D}}$ is strictly greater than one. Then, we have
$$
 \liminf_{n\to\infty} (F_n (u_n) )^{1/n} >1;
$$
hence, $F_n (u_n) \to \infty$ and the natural energy associated to
$u_n$ is unbounded.
\end{theorem}

\begin{proof}

Consider the absolutely minimizing Lipschitz extension (AMLE) of $f|_{\partial \Omega \cap \overline{D}
}$ to $D$, which is well defined even if the datum $f|_{\partial \Omega \cap \overline{D} }$ is not given in the whole $\partial D$. In this
case, the AMLE is characterized, as proved in \cite{PSSW} and
\cite{ChGAR}, as the unique solution of the problem
$$
\left\{\begin{array}{rcll}
- \Delta_{\infty} u (x)& = &0, \quad & x\in D,\\[5pt]
u(x) &= &f(x),\quad & x\in \partial \Omega \cap \overline{D},
\\[5pt]
\displaystyle \frac{\partial u}{\partial \nu}(x) & =&0, \quad & x
\in
\partial D \setminus \partial \Omega.
\end{array}\right.$$

Let $\lambda > 1$ be the Lispschitz constant in $D$ of this AMLE. Suppose that
$$
\liminf_{n \to \infty}  ( F_n(u_n) )^{1/n} = \beta < \lambda
$$
and consequently that
$$
\liminf_{n \to \infty}
\left( \int_D \frac{|\nabla u_n|^{n} }{n} \, dx \right)^{1/n} \leq
\beta.
$$
Fix $m \geq p_-$ and take $n > m$. By H\"{o}lder's inequality,
$$
\left( \int_{D} | \nabla u_n |^m \right)^{1/m} \le
|D|^{\frac{1}{m}-\frac{1}{n}} \left( \int_{D} | \nabla u_n |^n
\right)^{1/n} \le |D|^{\frac{1}{m}-\frac{1}{n}}n^{1/n} \left( \int_D
\frac{|\nabla u_n|^{n} }{n} \, dx \right)^{1/n}.
$$
Taking the limit in $n$, we conclude
$$
\liminf_{n\to \infty}\left( \int_{D} | \nabla u_n |^m \right)^{1/m} \le |D|^{\frac{1}{m}}
\beta,
$$
so, for a subsequence, there exists a weak limit in $W^{1,m}(D)$,
that we denote by $ u_{\infty} $.  This weak limit has to verify the
inequality
$$
\left( \int_{D} | \nabla u_{\infty} |^m \right)^{1/m} \le |D|^{\frac
1m} \beta
$$
for every $m$. Thus, taking the limit $m \to \infty$, we get that
$u_\infty \in W^{1,\infty} (D)$ and, moreover,
$$
| \nabla u_{\infty} | \le \beta, \qquad \textrm{a.e.}  \ \: x \in D.
$$
But this is a contradiction since $\lambda$ is the Lipschitz
constant in $D$ of the AMLE of $f|_{\partial \Omega \cap
\overline{D}}$ to $D$. We conclude that
$$
\liminf_{n \to \infty}  ( F_n(u_n) )^{1/n} \geq  \lambda
$$
and the result follows.
\end{proof}

\begin{remark}
The AMLE problem has been extensively studied in the
literature: see \cite{A}, \cite{J}, the survey \cite{ACJ}, and the
recent approach using tug-of-war games of \cite{ChGAR}, \cite{PSSW}
and \cite{PSSW2}.
\end{remark}

\begin{remark}\label{suf}
If $\partial \Omega \cap \overline{D} = \emptyset$, then $S \neq
\emptyset$; indeed, we can consider a function that is constant in
$D$ and coincides with $f$ on $\partial \Omega$, and extend it as a Lipschitz function to the whole of $\Omega$, thus obtaining an element of $S$.
\end{remark}

We now focus on the main case $S\neq \emptyset$. Recall that
solutions to \eqref{1.1}$_n$ are minima of the functional
$$F_n(u) = \int_\Omega \frac{|\nabla u|^{p_n (x)} }{p_n
(x)} \, dx
$$
in
$$
S_n = \left\{ u \in W^{1, p_n (x)} (\Omega) \, : \, u |_{\partial
\Omega}= f \right\}.
$$
The limit of these variational problems is given by minimizing
\begin{equation}\label{funcional.limite}
    F(u) =  \int_{\Omega \setminus \overline D} \frac{|\nabla u|^{p(x)} }{p(x)} \,
    dx
\end{equation}
in
$$
S = \left\{ u \in W^{1,p^{-}}(\Omega) \ : \ u|_{\Omega \setminus
\overline D} \in W^{1, p(x)} (\Omega \setminus \overline D), \ \ \|
\nabla u \|_{L^\infty (D)} \leq 1 \ \ \mbox{and} \ \ u |_{\partial
\Omega}= f \right\}.
$$

\begin{theorem} \label{teo.variacional}
Assume that $S\not= \emptyset$ and let $u_n$ be minimizers of $F_n$
in $S_n$. Then, along subsequences, $(u_n)_n$ converges uniformly in
$\overline{\Omega}$, weakly in $W^{1,m} (D)$, for every $m\geq p_-$,
and weakly in $W^{1,p(x)} (\Omega \setminus \overline D)$ to
$u_\infty$, a minimizer of $F$ in $S$. Moreover, the limit
$u_\infty$ is $\infty$-harmonic in $D$, \textit{i.e.},
$$
-\Delta_\infty u_\infty =0 \qquad \textrm{in} \ \: D,
$$
in the viscosity sense. Finally, the limit $u_\infty$ is unique, in
the sense that any other minimizer of $F$ in $S$ that is
$\infty$-harmonic in $D$ coincides with $u_\infty$.
\end{theorem}

\begin{proof}
We use the estimates obtained in the previous section. Since the
sequence $(u_n)_n$ is equicontinuous and uniformly bounded, by
Arzel\`a-Ascoli theorem it converges (along subsequences) uniformly
in $\overline{\Omega}$; the weak convergence in the space $W^{1,m}
(D)$, for every $m \geq p_-$, is obtained as in the proof of
Theorem~\ref{no.acotada} and the weak convergence in $W^{1,p(x)}
(\Omega \setminus \overline D)$ follows from the estimates in
Proposition \ref{lemma.estimaciones}.

Also as before, we get that $u_\infty \in W^{1,\infty} (D)$, with $|
\nabla u_{\infty} | \le 1$, a.e. $x \in D$, thus concluding that
$u_\infty \in S$. On the other hand, also from Proposition
\ref{lemma.estimaciones}, we get
$$ \int_{\Omega \setminus \overline D} \frac{|\nabla
u_n|^{p(x)}}{p(x)}  \, dx \leq  F_n(u_n) \leq F_n(v) \longrightarrow
\int_{\Omega \setminus \overline D}\frac{|\nabla v|^{p (x)} }{p (x)}
\, dx$$ and we conclude that
$$
F(u_\infty)=\int_{\Omega \setminus \overline D}\frac{|\nabla
u_\infty |^{p(x)} }{p (x)} \, dx \leq \int_{\Omega \setminus
\overline D}\frac{|\nabla v |^{p (x)} }{p (x)} \, dx =F(v), \quad
\forall \ v \in S
$$
so that $u_\infty$ is a minimizer for $F$ in $S$.

That a uniform limit of $n$-harmonic functions is $\infty$-harmonic
is a well known fact (\textit{cf.}, for example, \cite{BBM} or
\cite{J}).

To prove the uniqueness, suppose we have two minimizers in $S$,
$u_1$ and $u_2$. Then, considering
$$v = \frac{u_1+u_2}{2} \in S,$$ we obtain
that they coincide in $\Omega \setminus \overline{D}$ since $F$ is a
strictly convex functional in $S$. Using the uniqueness of solutions
of the Dirichlet problem for the $\infty$-Laplacian in $D$ (note
that $u_1$ coincides with $u_2$ on the whole of $\partial D$), we
conclude that $u_1=u_2$ also in $D$. We conclude, in particular,
that the whole sequence $u_n$ converges uniformly in
$\overline{\Omega}$.
\end{proof}

Our next task is to pass to the limit in \eqref{1.1.n.2}, the
problem satisfied by $u_n$ in the viscosity sense, to identify the
problem solved by $u_\infty$. We are under the assumption $S\neq
\emptyset$ and we recall that
$$u_n \to u_\infty$$
uniformly in $\overline{\Omega}$.

\begin{theorem} \label{teo.viscoso}
Every uniform limit of a sequence $\{u_n\}$ of solutions of
\eqref{1.1}$_n$ is a viscosity solution of
\begin{equation} \label{ec.viscosa}
\left\{\begin{array}{ll} -\Delta_{p(x)} u (x) =0, \qquad &
x\in\Omega
\setminus \overline D,\\
[5pt]-\Delta_{\infty} u (x) =0, \qquad & x\in D,\\
[5pt]
\displaystyle
 \mathrm{sgn} \left( |\nabla u | (x) -1\right) \: \mathrm{sgn} \left(
\frac{\partial u}{\partial \nu}(x) \right) =0, \qquad & x \in \partial D
\cap \Omega, \\[5pt]
u(x) = f(x),\quad & x\in \partial \Omega.
\end{array}\right.
\end{equation}
\end{theorem}

\begin{proof}
Since $u_n(x) = f(x)$, for $x \in \partial \Omega$, it is clear that
$ u(x) = f(x)$, for $x\in \partial \Omega$.

Let $u_{\infty}$ be a uniform limit of $\{u_n\}$ and let $\phi$ be a
test function such that $u_{\infty}(x_0)=\phi (x_0)$ and
$u_{\infty}-\phi$ has a strict minimum at $x_0\in \Omega$. Depending
on the location of the point $x_0$ we have different situations.

If $x_0 \in D$, we encounter the standard fact the the uniform limit
of $n$-harmonic functions is $\infty$-harmonic.

If $x_0 \in \Omega  \setminus \overline D$, consider  a sequence of
points $x_{n}$ such that $x_n \to x_0$ and $u_{n} -\phi$ has a
minimum at $x_{n}$, with $x_{n} \in \Omega  \setminus \overline D $
for $n$ large. Using the fact that $u_{n}$ is a viscosity solution
of \eqref{1.1.n.2}, we obtain
$$
- \Delta_{p_n(x_n)} \phi (x_{n}) \ge 0.
$$
Now we observe that $p_n(x) = p(x)$ in a neighborhood of $x_0$ and
hence, taking the limit as $n \to \infty$, we get
$$
- \Delta_{p(x_0)} \phi (x_{0}) \ge 0.
$$
That is, $u_{\infty }$ is a viscosity supersolution of
$-\Delta_{p(x)} u_{\infty} =0$ in $\Omega  \setminus \overline D $.

If $x_0 \in \partial D \cap \Omega$, we have to show that
$$\max \left\{
\begin{array}{l}
- \Delta_{p(x_0)} \phi (x_0) \\
[5pt] - \Delta_\infty \phi (x_0) \\
[5pt]   \mathrm{sgn} (|\nabla \phi | (x_0) -1) \: \mathrm{sgn} \left(
\frac{\partial \phi}{\partial \nu}(x_0) \right)
\end{array}
\right\}  \ge  0.
$$
Again, since $u_n$ converges to $u$ uniformly,  there exists a
sequence of points $x_n$ converging to $x_0$ such that $u_n - \phi$
has a minimum at $x_n$. We distinguish several cases.

\noindent \textsf{Case 1.} There exists infinitely many $n$ such
that $x_n \in D$.

\noindent Then we have, by
Proposition~\ref{lema.sol.debil.es.sol.viscosa},
$$
-\Delta_n \phi(x_n) =  - |\nabla \phi (x_n)|^{n-2} \Delta \phi
(x_{n}) -(n-2) |\nabla \phi(x_n)|^{n-4} \Delta_\infty \phi (x_{n})
\ge 0.
$$
If $\nabla \phi (x_0)=0$, we get $-\Delta_\infty \phi (x_0) =  0$.
If this is not the case, we have that $\nabla \phi (x_{n}) \neq 0$,
for large $n$, and then
$$
- \Delta_\infty \phi (x_{n}) \ge  \displaystyle\frac{1}{n-2} |\nabla
\phi (x_n)|^{2} \Delta \phi (x_{n}) \to 0, \quad \mbox{as } n \to
\infty .
$$
We conclude that
$$
-\Delta_\infty \phi (x_0) \geq 0.
$$

\noindent \textsf{Case 2.} There exists infinitely many $n$ such
that $x_n \in \Omega \setminus \overline D$.

\noindent Then we have, by
Proposition~\ref{lema.sol.debil.es.sol.viscosa},
$$
- \Delta_{p_n(x_n)} \phi (x_{n}) \ge 0.
$$
Proceeding as before, we get
$$
- \Delta_{p(x_0)} \phi (x_{0}) \ge 0.
$$

\noindent \textsf{Case 3.} There exists infinitely many $n$ such
that $x_n \in \partial D \cap \Omega$.

\noindent In this case, we have
$$
| \nabla \phi (x_n)|^{n-2}  \frac{\partial \phi}{\partial
\nu}(x_n)  - | \nabla \phi (x_n)|^{p(x_n)-2}  \frac{\partial \phi}{\partial \nu}(x_n)  \geq 0.
$$
Hence, we get
$$
\frac{\partial \phi}{\partial \nu}(x_n)  \leq  | \nabla \phi
(x_n)|^{n-p(x_n)} \frac{\partial \phi}{\partial \nu}(x_n).
$$
Taking $n\to \infty$, we deduce that
$$
|\nabla \phi | (x_0) > 1 \Rightarrow  \frac{\partial
\phi}{\partial \nu}(x_0) \geq 0,
$$
and
$$
|\nabla \phi | (x_0) < 1 \Rightarrow   \frac{\partial
\phi}{\partial \nu}(x_0)  \leq 0.
$$
That is
$$
\mathrm{sgn} (|\nabla \phi | (x_0) -1) \: \mathrm{sgn} \left(
\frac{\partial \phi}{\partial \nu}(x_0) \right) \geq 0.
$$
This concludes the proof that $u_\infty$ is a viscosity
supersolution.

The proof that $u$ is a viscosity subsolution runs as above and we
omit the details.
\end{proof}

\section{More on the set $S$}\label{remark222}
\setcounter{equation}{0}

We have already observed the following two facts concerning the non-emptiness of the set $S$:

\begin{enumerate}

\item If $\partial \Omega \cap \overline{D} = \emptyset$, then $S \neq
\emptyset$.

\medskip

\item If $\partial \Omega \cap \overline{D} \neq \emptyset$
and the Lipschitz constant of $f|_{\partial \Omega \cap
\overline{D} }$ is greater than one, then any Lipschitz extension
$u$ of this datum to $D$ verifies $\| \nabla u \|_{L^\infty (D)}
>1$ and, consequently, $S=\emptyset$.

\end{enumerate}

\noindent The question naturally arises of wether the
condition that the Lipschitz constant of  $f|_{\partial \Omega
\cap \overline{D}}$ is less than or equal to one is, not only
necessary, but also sufficient  to guarantee that $S \neq
\emptyset$.

Suppose we are given a Lipschitz boundary data $f$ such that the
Lipschitz constant of $f|_{\partial \Omega \cap
\overline{D} }$ is less than or equal to one. A natural attempt to
construct a function in $S$ would be the following:

\begin{itemize}

\item consider the unique AMLE of $f|_{\partial \Omega
\cap \overline{D}}$ to $D$,  which is such that the $L^\infty$-norm of its gradient
is less than or equal to one;

\medskip

\item extend it to the whole of $\Omega$ using any
function in $W^{1,p(x)} (\Omega\setminus \overline D)$ that
coincides with it on $\partial D$ and with $f$ on $\partial
\Omega$.
\end{itemize}

\noindent The boundary datum on $\partial
(\Omega \setminus D)$ that one has to extend is given by $f$
on $\partial \Omega \setminus \overline{D}$ and by the restriction
of the AMLE to $\partial D \cap \Omega$. The problem is that
the extension to $\Omega$ may not always be possible. However, if this
boundary data on $\partial (\Omega \setminus D)$ is Lipschitz,
then we could  indeed consider a Lipschitz extension to $\Omega$.

We first give an example of a particular geometric configuration
for which this is the case. Therefore the condition that the
Lipschitz constant of $f|_{\partial \Omega \cap
\overline{D} }$ is less than or equal to one does indeed suffice
to guarantee that $S \neq \emptyset$. Let $\Omega = B(0,1)$ in
$\RR^2$ and let
$$D= B(0,1) \cap \{ x > 0\}$$
be the right half-ball (here, $(x,y)$ denote coordinates in
$\RR^2$). We still denote by $f$ the obtained boundary datum on
the boundary of the half disc $B(0,1) \setminus D$, which is
Lipschitz on $\partial B(0,1) \cap \{x\leq 0\}$ and on $\partial D
= \{(x,y) \, :\, x =0, -1 \leq y \leq 1\}$, and continuous on the
whole boundary. Let $(0,y) \in \partial D \cap B(0,1)$ and $(z,w)
\in \partial B(0,1)$ with $w\geq 0$ and $y \geq 0$ (the other
possible cases would have to be considered separately). Adding and
subtracting $f((0,1))$ in the numerator we  obtain
$$
\frac{ |f((0,y)) - f ((z,w))| }{ \| (0,y) - (z,w) \| } \le \frac{ | f((0,y)) - f ((0,1)) | }{ \| (0,y) - (0,1) \| } + 2 \frac{ | f((z,w)) - f ((0,1)) | }{ \| (0,1) - (z,w) \| } \leq C,
$$
since
$\| (0,y) - (0,1) \| \leq \| (0,y) - (z,w) \| $, $\| (0,1) - (z,w) \| \leq  2 \| (0,y) - (z,w) \|$
and $f$ is Lipschitz on $\partial D$ and $\partial \Omega$. This shows that $f$ is Lipschitz on the whole boundary of $B(0,1) \setminus D$.

This construction does not always work since it may happen that the obtained boundary data is not a Lipschitz function. Here is a counter-example: let $\Omega = B((0,0),1)$ and $ D = B((1/2, 0), 1/2)$ in $\RR^2$. These two balls are tangent at the
point $(1,0)$. Now let $f$ be given in polar coordinates by
$$
f(\theta) = \left\{ \begin{array}{ll}
|\theta| \qquad & 0\leq |\theta| \leq \pi/2 \\
  \pi - |\theta| \qquad & \pi/2 <| \theta| < \pi.
\end{array}
\right.
$$
This function is Lipschitz on $\partial \Omega$. The unique AMLE
of $f|(1,0)$ to $D$ is given by $ u \equiv 0$. Now, we have the
function
$$
\tilde{f} (x) = \left\{ \begin{array}{ll}
f(x) \qquad & x\in \partial B((0,0),1) \\
0 \qquad & x \in \partial B((1/2,0), 1/2),
\end{array}
\right.
$$
defined on $\partial (\Omega \setminus D)$. Observe that there are points on $\partial D$ of the form $r (\theta) = \cos (\theta)$, with $\theta \sim 0$. For $\theta >0$,
$$
\lim_{\theta \to 0+}\frac{\tilde{f} (\cos (\theta), \sin (\theta)) - \tilde{f} (\cos
(\theta) (\cos (\theta), \sin(\theta)))}{ \|(\cos (\theta), \sin
(\theta)) - (\cos (\theta) (\cos (\theta), \sin(\theta)))  \|} =
\lim_{\theta \to 0+}
\frac{\theta-0}{ 1 - \cos (\theta)} = +\infty,
$$
hence the function is not Lipschitz.

\section{The one-dimensional case}\label{1-d}
\setcounter{equation}{0}

In this section, we analyze with some detail the one-dimensional
case, which is easier since the equation reduces to an ODE.

Let $\Omega = (0,1)$ and assume $p(x)\equiv \infty$ for $x\in (0,
\xi)$. Then the problem at level $n$ reads
$$
\left\{\begin{array}{l}
( | u_n' |^{p_n (x) -2 } u_n'  )' (x) =0, \\[8pt]
u_n(0) = f(0), \\[8pt]
u_n(1) = f(1).
\end{array}\right.
$$
To simplify, we assume that $f(0)=0$ and $f(1) >0$. Then,
integrating the equation, we get
$$
| u_n' |^{p_n (x) -2 } u_n' (x) = C_{1}.
$$
Assuming that $u_n' \geq 0$, we get
$$
 u_n' (x) = (C_1)^{\frac{1}{p_n (x) -1}}.
$$
Thus
$$
 u_n (x) = \int_0^x (C_1)^{\frac{1}{p_n (s) -1}}\, ds
$$
and the constant $C_1$ (that must be positive and depends on $n$)
verifies
$$
f (1) = \int_0^1 (C_1)^{\frac{1}{p_n (s) -1}}\, ds.
$$
Since $f(1)$ is finite, we conclude that $C_1$ must be bounded; if
not,
$$\lim_{n \to \infty} u_n (x) = u_\infty (x) =
+\infty$$ in the whole interval $(\xi,1]$ and this contradicts $u_n
(1) = f(1)$. Therefore, we can assume (taking a subsequence if
necessary) that
$$
\lim_{n \to \infty} C_1 (n) = C_\infty.
$$

\noindent \textsf{Case 1.} When $C_\infty >0$, we conclude that the
limit of $u_n$ is given by
$$
u_{\infty}(x)=\lim_{n \to \infty} u_n (x) = \left\{
\begin{array}{ll}
x \qquad & x \in [0,\xi], \\[10pt]
\displaystyle \xi + \int_\xi^{x} (C_\infty)^{\frac{1}{p(s) -1}}
\,ds\qquad & x \in [\xi, 1].
\end{array}
\right.
$$
As $u_n(1) = f(1)$, we realize that the constant $C_\infty$ is
determined by
$$
\xi + \int_\xi^{1} (C_\infty)^{\frac{1}{p(s) -1}} \,ds = f(1).
$$
This case, $C_\infty >0$, actually happens when $f(1) > \xi$. Since
$C_\infty$ is uniquely determined, we obtain the convergence of the
whole sequence $u_n$.

Note that in this case we can verify that $u_\infty$ is a minimizer
of the functional $F$ given by \eqref{funcional.limite}. Indeed,
since $|u'_{\infty}|(x) \leq 1$, for $x \in [0,\xi]$, we have that
$u_\infty \in S$ and since $u_\infty$ is a solution of
$$
( | u' |^{p (x) -2 } u'  )' (x) =0, \qquad u(\xi) = \xi,
\quad u(1) = f(1),
$$
we have that it minimizes the functional $F$, which in this case is
given by
$$
F(u_\infty) = \int_\xi^1 \frac{(C_\infty)^{\frac{p(s)}{p(s)
-1}}}{p(s)} \, ds,
$$
among functions that verify $u(\xi) = \xi$ and $ u(1) = f(1)$.

Now, for any function $w \in S$, we have $|w'|(x) \leq 1$, for $x
\in [0, \xi]$, and we get $w(\xi) \leq \xi$. Let $z$ be the solution
of
$$
( | z' |^{p (x) -2 } z'  )' (x) =0, \qquad z(\xi) = w(\xi)
\leq \xi,
\quad z(1) = f(1).
$$
Then we have
$$
F(w) \ge F(z) \geq F(u_\infty).
$$
To see that the last inequality is true just use the monotonicity
of the function
$$
C \mapsto \int_\xi^1 \frac{(C)^{\frac{p(s)}{p(s) -1}}}{p(s)} \, ds
$$
with respect to $C$.

\medskip

\noindent \textsf{Case 2.} When $C_\infty = 0$, we have that
$$
\lim_{n \to \infty} u_n (x) = \left\{
\begin{array}{ll}
K x \qquad & x \in [0,\xi], \\[10pt]
\displaystyle K \xi \qquad & x \in [\xi, 1].
\end{array}
\right.
$$
Here $K\leq 1$ is given by
$$
K = \lim_{n\to \infty} (C_1 (n) )^{1/n}
$$
(recall that we are taking $p_n(x) = p(x) \wedge n$).

As $u_n(1) = f(1)$ we get that the constant $K$ is given by
$$
K \xi = f(1).
$$
This case actually happens when $f(1) \leq \xi$. Since $K$ is
uniquely determined, we obtain the convergence of the whole sequence
$u_n$.

Note that in this case the limit $u_\infty$ is not differentiable,
but it is Lipschitz. Also note that it is easy to verify that
$u_\infty$ is a minimizer of the functional $F$ given by
\eqref{funcional.limite}. Indeed, $F(u_\infty) =0$ and $F(w) \geq
0$, for every $w\in S$.

\bigskip

\noindent {\bf Acknowledgements.} We thank Peter H\"ast\"o for a
conversation on variable exponent Sobolev spaces related to Lemma
2.1 and the referees for carefully reading the manuscript, and
making various suggestions that helped improving the exposition.

Part of this work was done during a visit of JDR to University of
Pittsburgh. He wants to thank for the friendly and stimulating
atmosphere found there.

JJM partially supported by NSF Award DMS-0500983. JDR partially
supported by project MTM2004-02223, MEC, Spain, by UBA X066 and by
CONICET, Argentina. JMU supported by CMUC/FCT.


\end{document}